\documentclass[10pt,conference]{IEEEtran}
\IEEEoverridecommandlockouts
\usepackage{graphicx}
\usepackage{amssymb,amsmath}
\usepackage{algorithm}
\usepackage{algpseudocode}
\usepackage{color}
\usepackage{xcolor} %Can remove. Only for comments. 
\usepackage{multirow} % To get two-lined rows in tables
\usepackage[normalem]{ulem}
\newcommand{\Input}{\textbf{Input: }}
\newcommand{\Output}{\textbf{Output: }}
\usepackage{hang} %To use hanging indent in algorithms
\newcommand{\R}{\mathbb{R}}
\DeclareMathOperator*{\argmin}{arg\,min}
\usepackage{physics} %Used in algorithms

\begin{document}

% special control bibtex entry
% won't be cited, but sets options for replacing long author lists with 'et al'
\bstctlcite{BSTcontrol}

%\title{Experimental Evaluation of Multiprecision Strategies for GMRES on GPUs\\
\title{Experimental Evaluation of Multiprecision Strategies for GMRES on GPUs\\
\thanks{%Center for Computing Research, Sandia National Laboratories, Albuquerque, NM, 87123 \{jloe,caglusa,iyamaza,egboman,srajama\}@sandia.gov
 Sandia National Laboratories is a multimission laboratory managed and operated by National Technology and Engineering Solutions of Sandia, LLC, a wholly owned subsidiary of Honeywell International, Inc., for the U.S. Department of Energy's National Nuclear Security Administration under contract DE-NA-0003525. This paper describes objective technical results and analysis. Any subjective views or opinions that might be expressed in the paper do not necessarily represent the views of the U.S. Department of Energy or the United States Government.  SAND2021-3861 C
 }
}

% \author{\IEEEauthorblockN{Jennifer Loe}
% \IEEEauthorblockA{\textit{Center for Computing Research} \\
% \textit{Sandia National Laboratories}\\
% Albuquerque, NM, USA \\
% jloe@sandia.gov}
% \and
% \IEEEauthorblockN{Christian Glusa}
% \IEEEauthorblockA{\textit{Center for Computing Research} \\
% \textit{Sandia National Laboratories}\\
% Albuquerque, NM, USA \\
% caglusa@sandia.gov}
% \and
% \IEEEauthorblockN{Ichitaro Yamazaki}
% \IEEEauthorblockA{\textit{Center for Computing Research} \\
% \textit{Sandia National Laboratories}\\
% Albuquerque, NM, USA \\
% iyamaza@sandia.gov}
% \and
% \IEEEauthorblockN{Erik G. Boman}
% \IEEEauthorblockA{\textit{Center for Computing Research} \\
% \textit{Sandia National Laboratories}\\
% Albuquerque, NM, USA \\
% egboman@sandia.gov}
% \and
% \IEEEauthorblockN{Siva Rajamanickam}
% \IEEEauthorblockA{\textit{Center for Computing Research} \\
% \textit{Sandia National Laboratories}\\
% Albuquerque, NM, USA \\
% srajama@sandia.gov}
% }

\author{\IEEEauthorblockN{Jennifer A.\ Loe%\IEEEauthorre fmark{1}
, Christian A.\ Glusa%\IEEEauthorrefmark{1}
, Ichitaro Yamazaki%\IEEEauthorrefmark{1}
, Erik G.\ Boman%\IEEEauthorrefmark{1}
,\ and Sivasankaran Rajamanickam%\IEEEauthorrefmark{1}
} 
\IEEEauthorblockA{%\IEEEauthorrefmark{1}
Center for Computing Research\\
Sandia National Laboratories \\
Albuquerque, New Mexico, USA 87123\\
\{jloe, caglusa, iyamaza, egboman, srajama\}@sandia.gov} }

%\title{Multiprecision Strategies for GMRES on GPUs}
% \author{Jennifer Loe, Christian Glusa, Ichitaro Yamazaki, Erik Boman, Siva Rajamanickam \footnotemark[1]}
% \date{\today}

%\footnotetext[1]{Center for Computing Research, Sandia National Laboratories, Albuquerque, NM, 87123 \{jloe,caglusa,iyamaza,egboman,srajama\}@sandia.gov
 %Sandia National Laboratories is a multimission laboratory managed and operated by National Technology and Engineering Solutions of Sandia, LLC, a wholly owned subsidiary of Honeywell International, Inc., for the U.S. Department of Energy's National Nuclear Security Administration under contract DE-NA-0003525. This paper describes objective technical results and analysis. Any subjective views or opinions that might be expressed in the paper do not necessarily represent the views of the U.S. Department of Energy or the United States Government.}

\maketitle

%TODO Remove these two lines before submitting!
%This forces page numbers so I don't have to count...
% \thispagestyle{plain}
% \pagestyle{plain}

\begin{abstract}
Support for lower precision computation is becoming more common in accelerator hardware due to lower power usage, reduced data movement and increased computational performance. However, computational science and engineering (CSE) problems require double precision accuracy in several domains. This conflict between hardware trends and application needs has resulted in a need for multiprecision strategies at the linear algebra algorithms level if we want to exploit the hardware to its full potential while meeting the accuracy requirements. In this paper, we focus on preconditioned sparse iterative linear solvers, a key kernel in several CSE applications. We present a study of multiprecision strategies for accelerating this kernel on GPUs. We seek the best methods for incorporating multiple precisions into the GMRES linear solver; these include iterative refinement and parallelizable preconditioners. Our work presents strategies to determine when multiprecision GMRES will be effective and to choose parameters for a multiprecision iterative refinement solver to achieve better performance. We use an implementation that is based on the Trilinos library and employs Kokkos Kernels for performance portability of linear algebra kernels. Performance results demonstrate the promise of multiprecision approaches and demonstrate even further improvements are possible by optimizing low-level kernels.
\end{abstract}


\section{Introduction}
In the current push towards exascale, modern supercomputers are increasingly relying on accelerator hardware for improved performance (with few exceptions). 
%In modern computing, algorithms are often memory-bound rather than compute-bound \CG{[Too general maybe?]}, so speedup depends on reducing data movement. 
These accelerators are starting to support and even rely on lower precision computations as their primary use case. This is due to lower power usage, reduced data movement with lower memory footprint requirements, and increased computational performance for lower precision computations. The emergence of machine learning accelerators, such as Cerebras, Sambanova, and Graphcore, which support only lower precision, increases the adoption of lower precision even further. In addition to increased efficiencies, most of these accelerators are being designed to address the needs of machine learning use cases in the industry that can tolerate 32-bit or even 16-bit computations. 

Using lower precision is starting to become important to realize the full potential of emerging hardware.
However, computational science and engineering (CSE) problems have a need for 64-bit computations. This level of accuracy is important because several of these simulations are used for high-consequence decision making. 
This conflict between the hardware trend and the application requirements has resulted in a renewed interest in multiprecision algorithms at the linear algebra library level \cite{AnztMPOverview}.
%where storing data in lower precision can reduce the memory load and take advantage of new hardware specifically designed for lower precisions, such a GPU tensor cores. 
Large-scale physics simulations with multiple discretized partial differential equations (PDEs) are also looking to take advantage of lower data precisions; however, unlike in machine learning, it is not obvious how to incorporate low precision data in the algorithm while obtaining double-precision accuracy of the final solution.

We focus on one of the expensive portions of solving PDEs, the sparse linear solve. While there are several approaches for solving sparse linear systems, we focus on sparse iterative linear solvers. The conjugate gradient (CG) method is highly effective for symmetric positive definite linear systems $Ax=b$. In this paper, we focus on the Generalized Minimum Residual method (GMRES) \cite{SaadSchultzGMRES}, which is commonly used for nonsymmetric systems. 

One algorithm that shows promise for this particular problem is GMRES with iterative refinement (GMRES-IR) \cite{WalkerTurnerGM-IR}. While the algorithm is several decades old, recent work with promising new analysis \cite{CarsonHighamNewAnalysis, CarsonIR3Precision} of this approach has increased interest. %\JL{Do we need to cite more refs here?}
%that has showed the promise of this approach with new analysis \cite{CarsonHighamNewAnalysis} has increased the interest in this approach.
However, this method has not been well-studied on modern accelerator-based architectures, and the algorithm is not standard in linear solver software implementations. % in frameworks in Trilinos. 
We address this gap by developing a Trilinos-based implementation of GMRES-IR. We further use this implementation for an experimental study that demonstrates the benefit of using GMRES-IR and, in some cases, what more needs to be improved. 

The main contributions of the paper are: 
\begin{itemize}
 \item Experimental evaluation of a Trilinos-based implementation of GMRES and two multiprecision variants, GMRES-FD and GMRES-IR, on GPUs; they show the promise of GMRES-IR for large problems that could take hundreds of iterations to converge.
 \item A demonstration with both model problems and general problems from the Suitesparse collection that GMRES-IR could reduce solve time by up to $1.5\times$ for preconditioned problems and $1.4\times$ for non-preconditioned problems while maintaining double precision accuracy. 
 \item An in-depth analysis of speedup of individual kernels within GMRES-IR on GPUs. 
 \item Evaluation of GMRES-IR combined with block Jacobi and polynomial preconditioning, and comparison with approaches such as low precision preconditioning with a higher precision solve.
 \item Evaluation of important GMRES-IR parameters such as subspace size and arithmetic complexity of preconditioners, as well as suggestions for tuning them for best performance.
\end{itemize}

 Our aim is that these experimental results will help users to have realistic expectations about potential performance gains from GMRES-IR, a starting place for parameter selection, and an understanding of effective preconditioning choices.

\section{Related Work}
%\JL{Come back and Reread this section.}
%\sr{This covers just two papers. Move other papers from intro here and cover them well. This doesn't have to even limited to GMRES-IR, you could also point to multiprecision work as well.}

The strategy of using low-precision computations to obtain high-precision solutions goes back (at least) to the 1960s. Recently, there has been a renewed interest in mixed-precision (multiprecision) methods~\cite{AnztMPOverview, BaboulinDongarraMPAlgs}. 
The most successful approach for linear systems has been \emph{iterative refinement}~\cite{Moler-IR}. The key idea is to compute $A \approx LU$ in low precision, which is both faster and requires less memory than the standard double precision factorization. Initially, one solves for $Ax^0 \approx LU x^0 = b$ in low precision, but then computes the residual $r^k =b-Ax^k$ in high precision and solves for a correction term using the error equation $A \Delta x^k = r^k$. By updating the previous solution by the correction term, $x^{k+1} = x^k + \Delta x^k$, a more accurate solution is obtained. One can iterate (reusing the $LU$ factors) until the desired accuracy is reached, typically in just a few iterations.
Iterative refinement has been highly successful for dense systems, especially on GPUs~\cite{HaidarMPIterRefGPU,DongarraGPUIR}.

We focus on iterative methods for sparse systems, which do not require $LU$ factorization. Several recent works have studied using multiple precisions with GMRES, including \cite{ CarsonHighamNewAnalysis,CarsonIR3Precision, HaidarMPIterRefGPU,DongarraGPUIR,AnztMPIterRef,  GrattonExploitingVPGMRES,  LindquistGMRES}. Anzt et al. \cite{AnztMPIterRef} analyzed iterative refinement combined with iterative solvers, viewing them as inner-outer solvers with iterative refinement as the outer solver and Krylov methods as the inner solver. They also presented some empirical results.

The original GMRES algorithm \cite{SaadSchultzGMRES} assumes every computation is done in high precision.
Turner and Walker \cite{WalkerTurnerGM-IR} observed that only a few key computations (including the residual) need to be done in high precision, while the rest can be done in lower precision. This approach has recently been revived as GMRES-IR \cite{CarsonHighamNewAnalysis,CarsonIR3Precision}.
A related approach is to compute the GMRES orthogonalization in lower (variable) precision~\cite{GrattonExploitingVPGMRES}.
Another option is \emph{inexact Krylov} methods~\cite{SzyldInexactKrylov}, but this was designed for inexact matrix-vector products and it is difficult to adapt to our mixed (single, double) precision use case.

The typical GMRES-IR implementations studied by Carson and Higham \cite{CarsonHighamNewAnalysis,CarsonIR3Precision} used various $LU$ factorizations in low precision as a preconditioner. There are two drawbacks of this approach. First, exact $LU$ may require too much memory due to fill (in the sparse case), has high computational complexity (due to fill), and may not be practical for large systems. Second, these preconditioners require a global triangular solve, which is not highly parallelizable, so not suitable for GPUs. Therefore, we do not consider $LU$-types of preconditioning here. 
The experiments in recent work were limited to small problems in MATLAB on CPUs. 
Instead, we focus on classical sparse preconditioners such as block Jacobi and matrix polynomials, which are more efficient on GPUs. 
\emph{The key contribution of this work is to evaluate this algorithm that shows promise in theory on a hardware that is designed to do well when using lower precision computation.} 
%In other words our evaluations of GMRES-IR is more general as it allows any preconditioner in GMRES, not just $LU$. 
The GMRES-IR algorithm we consider is given in Algorithm \ref{alg:GMRES-IR} and is essentially the method by Turner and Walker \cite{WalkerTurnerGM-IR}.

Recently, in concurrent work, an empirical study of GMRES and GMRES-IR with incomplete LU factorizations \cite{LindquistGMRES} was presented. While it is similar in scope, the evaluations there are limited to CPUs. In addition to an experimental study on GPUs, we also use different preconditioners, provide kernel-level performance analysis, and compare GMRES-IR to other schemes such as a precision switching scheme.

\section{GMRES and Multiprecision Variants}
% MGS Version:
% \begin{algorithm}[t]
%  \caption{GMRES(m) (MGS) \cite[p.\ 172]{SaadItMeth} }
%  \begin{hangingpar}
%  \Input $A\in \R^{n\times n}$, $b \in \R^{n\times 1}$, initial guess $x_0\in R^{n\times 1}$, relative residual tolerance $rTol$
%  \end{hangingpar}
%   \Output approximate solution $x_m$
%   %\hspace*{\algorithmicindent}
%   \begin{algorithmic}[1]
%     \State $r_0 = b-Ax_0$,
%     \State $\gamma = \|{r_0}\|_2$ and $v_1 = r_0/\gamma$ 
%       \For{$j=1:m$}
%       \State $w_j = Av_j$
%         \For{$i=1:j$}
%           \State $h_{ij} = v_i^T w_j$
%           \State $w_j = w_j - h_{ij}v_i$
%         \EndFor
%       \State $h_{j+1,j} = \|{w_j}\|_2$. (Lucky breakdown if $h_{j+1,j}=0$.)
%       \State $v_{j+1}= w_j/h_{j+1,j}$
%       \EndFor
%     \State Define the $(m+1)\times m$ upper-Hessenberg matrix $\overline{H}_m = \{h_{ij}\}_{1\le i\le m+1, 1\le j \le m}$
%     \State Compute $\hat{d} = \argmin_{y\in \R^{m}} \|{\gamma e_1 - \overline{H}_m y}\|_2$, $\hat{x} = V_m\hat{d}$, and $x_m = x_0 + \hat{x}$.
%     \State Compute $r_m = b - Ax_m$. If $\|{r_m}\|_2/\|{r_0}\|_2 \le rTol$, stop. Else, set $x_0 = x_m$, $r_0 = r_m$ and go to Step 2. 
%   \end{algorithmic}
%   \label{alg:GMRES(m)}
% \end{algorithm}

%FYI- If I needed to write CGS with loop, it would be:
%       \State $w_j = Av_j$
%       \State $tmp = w_j$
%         \For{$i=1:j$}
%           \State $h_{ij} = v_i^T tmp$
%           \State $w_j = w_j - h_{ij}v_i$
%         \EndFor
% Don't need a tmp in blocked version because wj isn't changed
% until all dot prods are done. 

%CGS version with blocking:
\begin{algorithm}[t]
  \caption{GMRES(m) (CGS) \cite[p.\ 172]{SaadItMeth} }
  \begin{hangingpar}
  \Input $A\in \R^{n\times n}$, $b \in \R^{n\times 1}$, initial guess $x_0\in R^{n\times 1}$, relative residual tolerance $rTol$
  \end{hangingpar}
  \Output approximate solution $x_m$
  %\hspace*{\algorithmicindent}
  \begin{algorithmic}[1]
    \State $r_0 = b-Ax_0$,
    \State $\gamma = \|{r_0}\|_2$, $v_1 = r_0/\gamma$, and $h_{1,1} = 0$
    \State Let $H_{:,j}$ be the vector of elements $\{h_{i,j}\}_{1\le i \le j}$.
      \For{$j=1:m$}
      \State $w_j = Av_j$
      \State Define $V_j = [v_1,v_2, \ldots, v_j]$. 
      \vspace{0.02in}
      \State $H_{:,j} = V_j^T w_j$ %(\texttt{GEMV Trans})
      \State $w_j = w_j - V_j H_{;,j}$ %(\texttt{GEMV No Trans})
      \State $h_{j+1,j} = \|{w_j}\|_2$. (Lucky breakdown if $h_{j+1,j}=0$.)
      \State $v_{j+1}= w_j/h_{j+1,j}$
      \EndFor
    \State Define matrix $\overline{H}_m = \{h_{i,j}\}_{1\le i\le m+1, 1\le j \le m}$.
    \State Compute $\hat{d} = \argmin_{y\in \R^{m}} \|{\gamma e_1 - \overline{H}_m y}\|_2$, $\hat{x} = V_m\hat{d}$, and $x_m = x_0 + \hat{x}$.
    \State Compute $r_m = b - Ax_m$. If $\|{r_m}\|_2/\|{r_0}\|_2 \le rTol$, stop. Else, set $x_0 = x_m$, $r_0 = r_m$ and go to Step 2. 
  \end{algorithmic}
  \label{alg:GMRES(m)}
\end{algorithm}

We begin by describing GMRES, the computational kernels involved, and important observations for a multiprecision approach. We follow this with a description of two multiprecision variants, GMRES-IR and GMRES-FD.

\subsection{GMRES}
%\sr{Point out line numbers on what you are describing}
We consider real-valued $n\times n$ sparse linear systems $Ax=b$. GMRES(m) (Algorithm \ref{alg:GMRES(m)}) builds out a Krylov subspace $\mathcal{K}_m(A,b)=\text{span}\{b, Ab, A^2b, \ldots, A^{m-1}b\}$ from which to extract an approximate solution $\hat{x}$. % for the system $Ax=b$. 
At each iteration, GMRES appends a new basis vector to the subspace, orthogonalizes that vector against the previous basis vectors, and uses the expanded subspace to update the approximate solution $\hat{x}$. 

GMRES has ``converged" when the relative residual norm $\|b-A\hat{x}\|_2/\|b\|_2$ falls below some user-specified tolerance. We say that GMRES convergence has improved when either (a) the total solve time decreases or (b) the iteration count for convergence decreases. When computing with only one precision, (a) and (b) are roughly equivalent, but this will not always be the case when comparing double precision (fp64) GMRES with a multiprecision implementation. 

GMRES is optimal in the sense that it picks the approximate solution $\hat{x}$ so that the residual norm $\|b-A\hat{x}\|_2$ is minimized. When the dimension of the Krylov subspace becomes too large (i.e.\ orthogonalizing a new basis vector becomes too expensive or the set of $m$ basis vectors of length $n$ can no longer fit in memory), we \emph{restart} GMRES. This means that we discard the current Krylov subspace and start the GMRES iteration from the beginning with the new right-hand side $r=b-A(x_0-\hat{x})$. Then the final solution is the sum of the initial starting vector $x_0$ and all intermediate solution vectors $\hat{x}$. We refer to the value $m$ as the \textit{maximum subspace size} or the \textit{restart length} for GMRES.

Note that restarting GMRES can slow convergence. When restarted, GMRES loses crucial eigenvector information from the previous subspace that allows it to converge more quickly to a solution \cite{MorganGMRESE}. It has to recreate this information in the next subspace, which requires more time and iterations. Thus, it can be a challenge to choose a restart length for GMRES that is large enough for quick convergence but small enough to fit in memory on GPU accelerators. 

%\paragraph{Computation Kernels} Don't think this makes sense anymore -J
The primary sources of computational expense for GMRES are (1) sparse matrix-vector products (SpMVs) with the matrix $A$ (Alg.\ \ref{alg:GMRES(m)}, line $5$) and (2) orthogonalization of the Krylov subspace vectors. In our experiments, each GMRES iteration uses two passes of classical Gram-Schmidt orthogonalization (CGS2). Each of these two orthogonalization passes requires two calls to GEMV, one with a transpose to compute inner products, and another with no transpose to subtract out components of the previous vectors (Alg.\ \ref{alg:GMRES(m)}, lines $7$ and $8$). Other less expensive operations include norms, small dense matrix operations with matrix $H$, and vector additions. 

%\paragraph{Notes}

\subsection{Multiprecision GMRES-IR}
For GMRES with iterative refinement (GMRES-IR), we will run the GMRES algorithm in \emph{single precision} (fp32) and then ``refine" the algorithm at each restart by starting the next GMRES run with a right-hand-side vector that has been computed in \emph{double precision} (fp64). (See Algorithm \ref{alg:GMRES-IR}.) Thus, we maintain both double and single precision copies of the matrix $A$ in memory for performing SpMVs in the appropriate precision. 
Note that we only check for convergence of GMRES-IR at each restart, when the residuals are recomputed. This is different from standard GMRES where we can monitor an implicit residual within the iteration to alert us to convergence. For GMRES-IR, the single precision residuals of the inner GMRES solver give little information about the convergence of the overall problem. Thus, GMRES-IR may take at most $m-1$ extra iterations in single precision over what is absolutely needed for convergence. We include cost of such iterations in our performance comparisons.
\begin{algorithm}
\caption{GMRES-IR}\label{alg:GMRES-IR}
  \begin{algorithmic}[1]
  \State $r_0 = b-Ax_0$ [double]
  \For{ $i=1,2, \ldots$ until convergence:}
    \State Use GMRES$(m)$ to solve $Au_i = r_i$ for correction $u_i$ [single]
    \State $x_{i+1} = x_i + u_i$ [double]
    \State $r_{i+1} = b - Ax_{i+1}$ [double]
  \EndFor
  \end{algorithmic}
\end{algorithm}

\subsection{Multiprecision GMRES-FD}

A first inclination when attempting to incorporate low precision into GMRES$(m)$ is to perform the entire first part of the calculation in one precision and then switch precisions at one of the restarts. We briefly explore these possibilities and then demonstrate why GMRES-IR is the better candidate for incorporating low precision. 

There are two options for switching precisions mid-solve: (1) Start in single precision and later switch to double, or (2) start in double and switch to single precision. The theory of \emph{inexact Krylov} supports option (2), stating that one can loosen the accuracy of the matrix-vector multiply (SpMV) as the iteration progresses and still converge to the correct solution \cite{SzyldInexactKrylov}. Furthermore, \cite{GrattonExploitingVPGMRES} shows that one can loosen the accuracy of inner products in addition to accuracy of the SpMV and still get convergence behavior close to that of full double precision. 
Option (2) may also be preferable because the initial computations in double precision may allow the Krylov subspace to quickly get good approximations to key eigenvectors, which can aid convergence \cite{MorganGMRESE}. However, inexact Krylov theory assumes that the vector operations are done in full precision and only the matrix-vector multiply is inexact. Therefore, the theory does not cover the use case of switching from double to single precision. 
%\JL{Not sure about the previous two sentences because we just said we could loosen accuracy of the inner products??} 
It is not clear if a single precision solver can even converge to double precision accuracy; thus, we do not evaluate option (2) in our experiments. We assess option (1), switching from single precision GMRES$(50)$ to double precision GMRES$(50)$ and using the single precision solution vector as a starting vector for the double precision GMRES iteration. We call this method GMRES-FD (Float-Double).

\subsection{Preconditioning} We investigate two lower precision alternatives to traditional fp64 preconditioning: (a) double precision GMRES with a single precision preconditioner and (b) GMRES-IR with a single precision preconditioner. Most previous studies of GMRES-IR (e.g.\ \cite{CarsonIR3Precision}) used some variation on LU preconditioning. %The LU factors were stored in lower precision so that the factorization was not exact \cite{CarsonIR3Precision}.
%Because large triangular solves are highly sequential, we do not believe this method is a good first choice for GPUs.
Here, we investigate more paralellizable preconditioners, using a polynomial preconditioner (Sections \ref{subsec:precon_compar} and \ref{subsec:poly}) and a block Jacobi preconditioner (Section \ref{subsec:suitesparse}). % for general problems.
In all tests, we use right preconditioning ($AMM^{-1}x = b$) so that the residuals of the preconditioned problem match those of the unpreconditioned problem in exact arithmetic. Each time an fp32 preconditioner $M$ is applied to an fp64 vector $x$ in case (a), we must cast $x$ to fp32, multiply it by $M$ in fp32, and cast the result back to fp64. For case (b), $M$ is both computed and applied entirely in fp32. 
 
 Polynomial preconditioning is applied as follows: We use a polynomial preconditioner based upon the GMRES polynomial (see details in \cite{LoePPTrilinos}). Here, using a polynomial $p(t)=\sum_{k=0}^{d}c_kt^k$ of degree $d$ as a preconditioner $M$ is to be understood as $M = p(A) = \sum_{k=0}^{d}c_k A^{k}$. (See \cite{XiaoPPCG} for a related study with the conjugate gradient method run in double precision with a single precision polynomial preconditioner.) 

\section{Software Implementation}

Trilinos \cite{Trilinos} is a large software library with packages for PDE discretizations, linear and non-linear solvers, preconditioners, partitioners, and distributed linear algebra. We use the Trilinos framework for our solver implementation, with the eventual goal of making GMRES-IR available in the public codebase. Thus, we test GMRES and GMRES-IR within the framework of Belos \cite{ThornquistBelosAmesos2}, the Trilinos sparse iterative linear solvers package. Our final software version will be available with Tpetra-based linear algebra and will run with MPI over many CPUs and GPUs. For this paper, however, we only consider solvers on a single CPU/GPU. For the solvers' linear algebra backend, we use the Kokkos \cite{TrottKokkos} and Kokkos Kernels libraries, which provide portable, optimized linear algebra operations for GPUs. 

% Skip MPI notes:
%That is, we wanted to avoid using MPI because MPI might make implicit type conversions during some of its operations. We wanted to avoid this and be certain of what scalar types we were using. 
% \CG{Is it true that MPI will perform implicit type conversions? Aren't we sending everythings as MPI\_BYTE or somethings like that anyway?}

 The Belos linear solvers package does not contain its own implementation of linear algebra, but instead relies on abstracted linear algebra interfaces through the \texttt{Belos::MultiVectorTraits}. We created a Kokkos-based adapter for Belos, letting the Kokkos adapter inherit from \texttt{Belos::MultiVector}. All of the length $n$ basis vectors for the Krylov subspace are stored in \texttt{Kokkos::Views} and operated on via the MultiVector interface. The interface implements all the needed capabilities to solve linear systems $Ax=b$ with a single right-hand side. 
Belos' solvers are all templated upon a user-specified scalar type, so they can be run in either float or double precision. Thus, at first glance, it seems that they would be well-suited for multiprecision computations. However, these templates assume that all operations are carried out in the same scalar type; there are no current capabilities to mix and match precisions within a solver. In spite of this, it is possible to perform operations outside of a solver using a different precision from the one the solver uses. We do this in our GMRES-IR implementation: The code initializes a Belos GMRES solver in fp32. At each restart, we retrieve the current solution vector from the Belos solver and convert it to fp64. Then we compute the current residual, convert that residual vector back to fp32, and feed that residual to the fp32 GMRES solver as the next right-hand side. 

\textit{Limitations of current implementation:}
Since we use the existing Belos interface, any mixed precision operations that are internal to the solver must be handled entirely in the linear algebra adapter. 
In order to avoid this difficulty, we do not study variations of GMRES where internal kernels use lower precision, e.g.\ GMRES with mixed precision orthogonalization or low precision SpMVs. 
Additionally, the Belos linear solvers package was not designed with GPUs or other accelerators in mind: Belos requires that results of some GPU operations be stored in a dense matrix representation on host (\texttt{Teuchos::SerialDenseMatrix}). %, which can only exist on the CPU. 
This requires data movement between the GPU and CPU along with memory allocations that otherwise might be unnecessary. 
%\CG{I guess the data movement is fine, since you want to compute inner products on device, but use them on host. But it also requires constant reallocation on device.} 
Furthermore, the structure in Belos forces separate kernel launches for each GPU operation, while in a Kokkos-only implementation some of these operations could be fused. We plan to improve upon these limitations in future software upgrades of the Belos package.

\section{Experimental Results}
\label{sec:experiments}

% Note: To re-create these results, use the KokkosBelos branch, tag JMPResults

All experiments that follow are run on a node equipped with a Power 9 CPU that has 318 GB DDR3 RAM and a Tesla V$100$ GPU with $16$ GB GDDR5 RAM. %All on Weaver so far. 
We used GCC 7.2.0, CUDA 9.2.88, Kokkos and Kokkos Kernels 3.2.0 and Trilinos 13.1. 
%TODO: Is Kokkos running with OpenMP on the host end? Nope. From what I'm seeing in using top on test runs, I'm only using one CPU core. TODO double-check this. 
%CPU operations are run in serial on a single process. %I'm not sure that is true. I think the default has OpenMP, even though it looked like it was only using one proc?
All PDE test problems used in sections \ref{subsec:IRvsFD} to \ref{subsec:poly} were generated with finite difference stencils via the Trilinos Galeri package.

%That is, we perform $m$ iterations of GMRES in float precision, compute the new residuals in double precision, and then restart the GMRES iteration using the truncated residual in float precision. %Now this is stated in earlier sections. 
The following experiments are run as follows: Unless otherwise stated, we restart both double precision GMRES$(m)$ and GMRES-IR after each run of $m=50$ iterations. 
%This info is earlier: 
%(For GMRES-IR, we recompute the residuals at each restart to check for convergence and apply the double precision refinement.) 
%Each GMRES iteration uses two passes of classical Gram-Schmidt (CGS2) orthogonalization. This allows orthogonalization operations to be grouped into four calls to GEMV. 
%Furthermore, CGS2 orthogonalization is more stable in maintaining an orthogonal basis at low precision. (TODO Cite?) %TODO Delete previous sentence if you can't cite it....
%In our experiments, GMRES will always restart after each $m=50$ iterations unless otherwise noted. 
%For our experiments we say that GMRES has converged when the relative residual norm drops below $1\mathrm{e}{-}10$. 
%Furthermore, we begin all our experiments with starting vector $x_0$ as all zeros. 
All solvers are run to a relative residual convergence tolerance of $1\mathrm{e}{-}10$.
For each problem, we use a right-hand side vector $b$ of all ones and a starting vector $x_0$ of all zeros. For each set of results, we exemplify the run that has the median of three solve times.
For GMRES-IR , total solve times do not include the time needed to make a single precision copy of matrix $A$, but they do include time required to convert residual vectors from double to single precision (and vice-versa) during the refinement stage. Note that results are not entirely deterministic; numerical errors from reductions on the GPU can give slightly different convergence behaviors. 

The rest of the experiments section is organized as follows. We compare different approaches for multiprecision GMRES (\ref{subsec:IRvsFD}). We evaluate GMRES and multiprecision GMRES-IR unpreconditioned (\ref{subsec:unprecon_compar}) and preconditioned (\ref{subsec:precon_compar}) for their convergence and performance. Next we do an in-depth analysis of the performance we observe in SpMV (\ref{sec:SpMV}). We also study how to choose two important parameters that affect performance: restart size (\ref{subsec:restart}) and arithmetic complexity of preconditioning (\ref{subsec:poly}). Finally, we evaluate our approach on a few general problems from the SuiteSparse collection (\ref{subsec:suitesparse}). 

\subsection{GMRES-IR vs GMRES-FD}
\label{subsec:IRvsFD}

\begin{figure}
  \centering
  \includegraphics[width=0.4\textwidth,keepaspectratio]{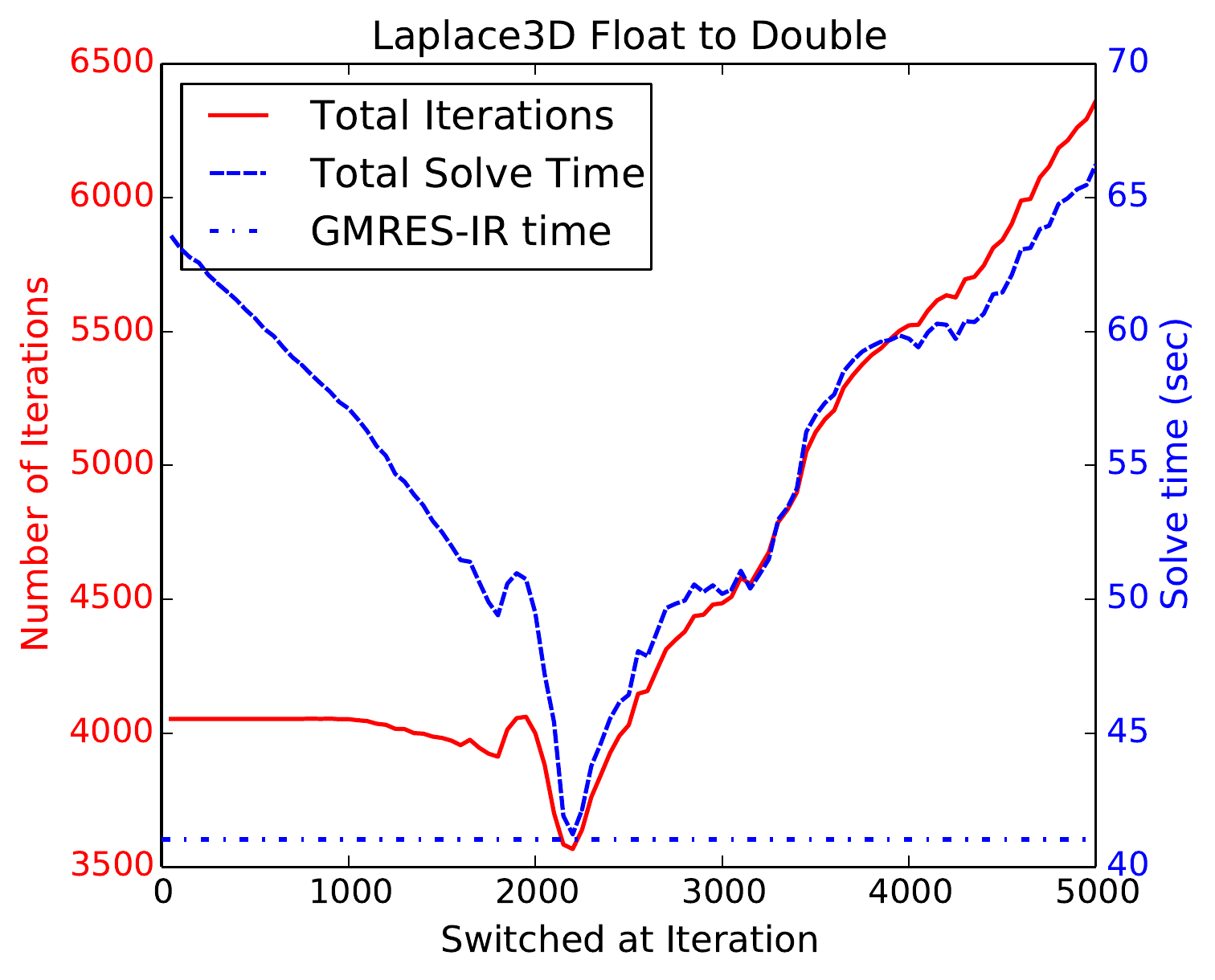}
  \caption{Total solve time and iteration count for a 3D Laplacian with GMRES-FD, switching from single precision to double precision at the iteration indicated on the horizontal axis. Dotted line at bottom indicates solve time for GMRES(50)-IR.}
  \label{fig:LaplSwitchIter}
\end{figure}

\begin{figure}
  \centering
  \includegraphics[width=0.4\textwidth,keepaspectratio]{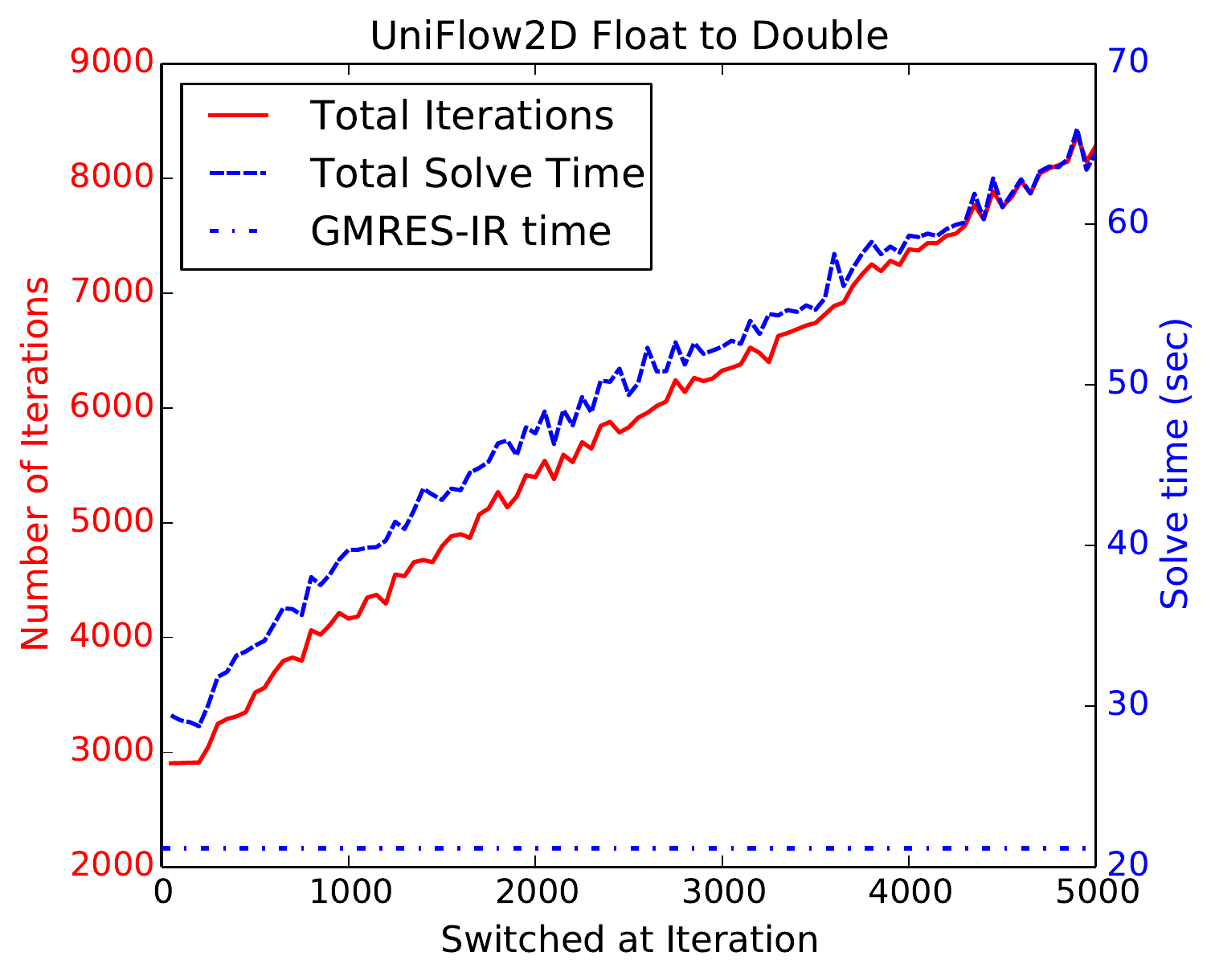}
  \caption{Total solve time and iteration count for the problem UniFlow2D2500 with GMRES-FD, switching from fp32 to fp64 at the iteration indicated on the horizontal axis. Dotted line at bottom indicates solve time for GMRES(50)-IR.}
  \label{fig:UniFlowSwIter}
\end{figure}

% But since a single precision solver will typically not converge to double precision accuracy, ending the computation in single precision is not helpful. %TODO triple check that this is true.
We begin experimental evaluations by comparing GMRES(m) in double precision, GMRES-IR, and GMRES-FD.
The first question with GMRES-FD: At what point is the right moment to switch precisions? We investigate with two different problems, comparing \emph{multiple runs of GMRES-FD} (switching at different iteration numbers), with a single run of GMRES-IR and GMRES(m). The first problem is a Laplacian from a 3D finite difference stencil with grid size $200$, and the second is a 2D convection-diffusion problem named ``UniFlow" with grid size $2500$. For both of these problems, we tested GMRES-FD, switching from fp32 to fp64 at each multiple of $50$ iterations (so at each restart). The $x$-axis in Figures \ref{fig:LaplSwitchIter} and \ref{fig:UniFlowSwIter} indicates the iteration at which the solver switched from float to double precision. The left vertical axis gives the total number of iterations required for convergence (the sum of single and double precision iterations). The right vertical axis gives total solve time for the problem. 
%Note that these results were not entirely deterministic; numerical errors from reductions on the GPU likely contributed to slightly different solution vectors from the single precision run and different starting vectors can greatly affect the convergence of GMRES double. %TODO: Do we actually need this line? 

One can predict that switching to fp64 too early is not harmful to convergence, but it does not take full advantage of the fp32 solver to find the minimum solve time. 
%The double precision solver may have to make extra iterations to make up for the lack of progress in single precision. % Trying to explain here why the iterations go down at the beginning- don't think this line quite does it. 
If the chosen switching point is too late, then the fp32 solver takes extra iterations, adding to the total solve time but not making any progress. This is exactly what we see with the Laplacian problem in Figure \ref{fig:LaplSwitchIter}. The solve time slowly for GMRES-FD decreases until reaching a minimum when the switch happens at $2200$ iterations. Here, the total number of iterations required for the solve is $3567$, while the solve time is $41.22$ seconds. 
%\sr{Can you add a x in Fig 1 and 2 for double precision GMRES and GMRES-IR number of iters and time?} 
Comparatively, GMRES$(50)$-IR converges in $4100$
iterations and $41.03$ seconds. The double precision-only problem requires $4053$ iterations and $63.83$ seconds. Thus, GMRES-IR attains the minimum solve time of all methods without needing to manually determine when to switch precisions. %Results from UniLaplIRTest1_14.out %TODO Would be nice to have a 'why' here: WHY would we expect GMRES-IR to find that minimum? Any intuition? 

Results from testing various switching points for GMRES-FD on the UniFlow problem (Figure \ref{fig:UniFlowSwIter}) are somewhat counterintuitive. The minimum of $28.77$ seconds (with a total of $2911$ iterations) occurs when switching at only $200$ iterations. This gives little improvement over the purely double precision solver, which required $2905$ iterations and $29.62$ seconds. Did the single precision solver's convergence stall after only $200$ iterations? Not at all! At a switching point of $2800$ iterations, for instance, the initial vector $x_0$ from the fp32 solver helps the fp64 solver to start with an initial residual norm of $9.9\mathrm{e}{-}5$. However, even with the good starting vector, the fp64 solver still needs an additional $3295$ iterations to converge. We hypothesize that this is because the new $x_0$ used at the switch of precisions did not contain eigenvector components that were present in the original right-hand side $b$. 
GMRES-IR, on the other hand, converges in $3000$ iterations and only $21.17$ seconds. It is the best method by far. This experiment demonstrates a case where GMRES-IR is quite helpful and GMRES-FD is mostly ineffective. %Results from UniLaplIRTest1_14.out 
We will use GMRES-IR as the multiprecision approach for the rest of the paper.
Next we look at how convergence of GMRES-IR compares to GMRES double and which kernels contribute most to speedup. 

%TODO: Unify capitalization of "Figure", "Table", etc. 

%\textbf{TODO: when things go wrong?}
%Wasn't there a matrix that was so ill-conditioned that GMRES just died in single precision? But I thought it was a bug? \CG{Biharmonic equation?} \JL{No, not Biharmonic. Xenon2 is one example, but there was also another.} %TODO investigate further. 

\subsection{Convergence and Kernel Speedup for GMRES vs GMRES-IR}
\label{subsec:unprecon_compar}
\begin{figure}
  \centering
  \includegraphics[width=0.45\textwidth,keepaspectratio]{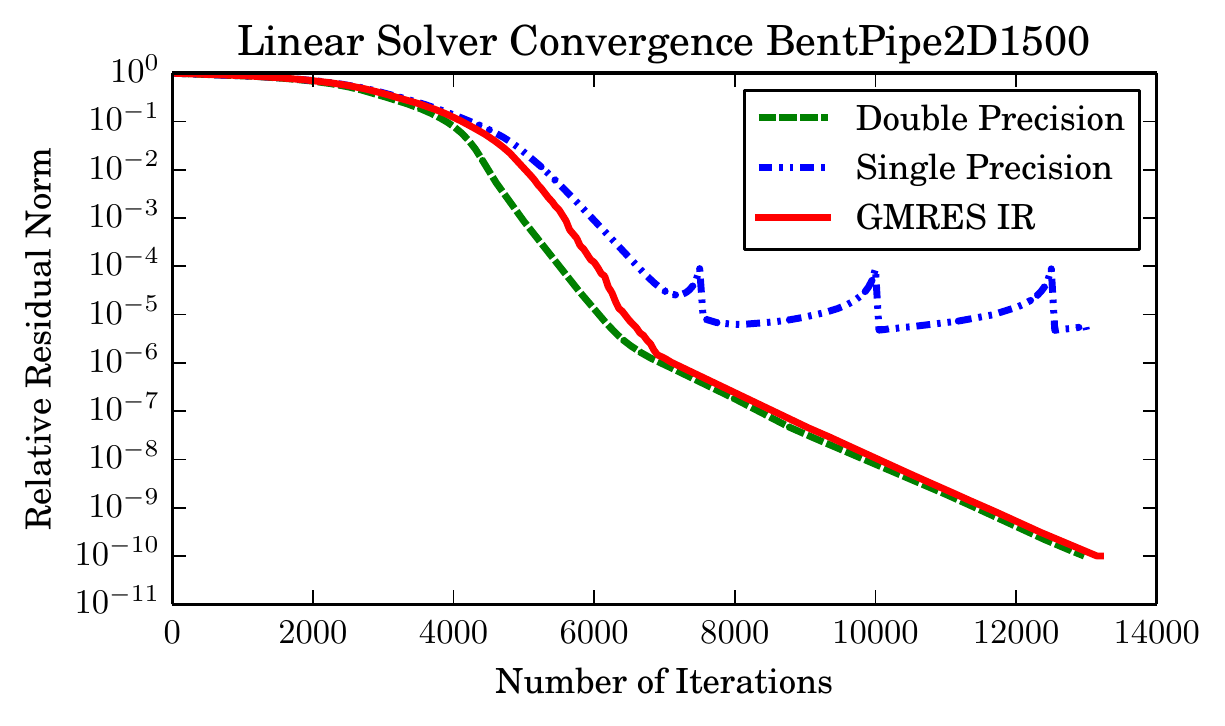}
  \caption{Relative residual norm convergence for matrix BentPipe2D1500. Single precision GMRES$(50)$ is represented by the blue dash-dot line, double precision by green dashes, and mixed precision GMRES$(50)$-IR by the red solid line.}
%\caption{Relative residual norm convergence for the matrix BentPipe2D1500. Single precision GMRES$(50)$ is represented by blue circles, double precision by green squares, and mixed precision GMRES$(50)$-IR by red triangles.}
  \label{fig:BentPipeConv}
\end{figure}

We next consider BentPipe2D1500, a 2D convection-diffusion problem with $nx=1500$, $n=2{,}250{,}000$ and $nnz = 11{,}244{,}000$. (Here $nx$ denotes the number of grid points in each direction of the mesh for the finite difference discretization of the PDE, and $nnz$ denotes the number of nonzero elements in the sparse matrix $A$.) 
The underlying PDE is strongly convection-dominated, so the matrix is ill-conditioned and highly non-symmetric.
We compare GMRES$(50)$ in all single precision, GMRES$(50)$ in all double precision, and GMRES$(50)$-IR. Convergence plots are in Figure \ref{fig:BentPipeConv}. 
The fp32 solver reaches a minimum relative residual norm of about $4.7\mathrm{e}{-}6$, and the fp64 solver needs $12{,}967$
iterations to converge to $1\mathrm{e}{-}10$. 
GMRES-IR needs $263$ cycles of $50$ iterations to converge, so $13{,}150$ total iterations, and its convergence curve closely follows that of the double precision solver. This phenomenon is related to the theory built by \cite{GrattonExploitingVPGMRES} for non-restarted GMRES; it has also been observed by \cite{LindquistGMRES} for restarted GMRES.
To reiterate, \emph{the convergence of the multiprecision version of the solver follows the double precision version closely}. 
% \textbf{GMRES-IR convergence tracks well with double precision}

\begin{figure}
  \centering
  \includegraphics[width=0.45\textwidth,keepaspectratio]{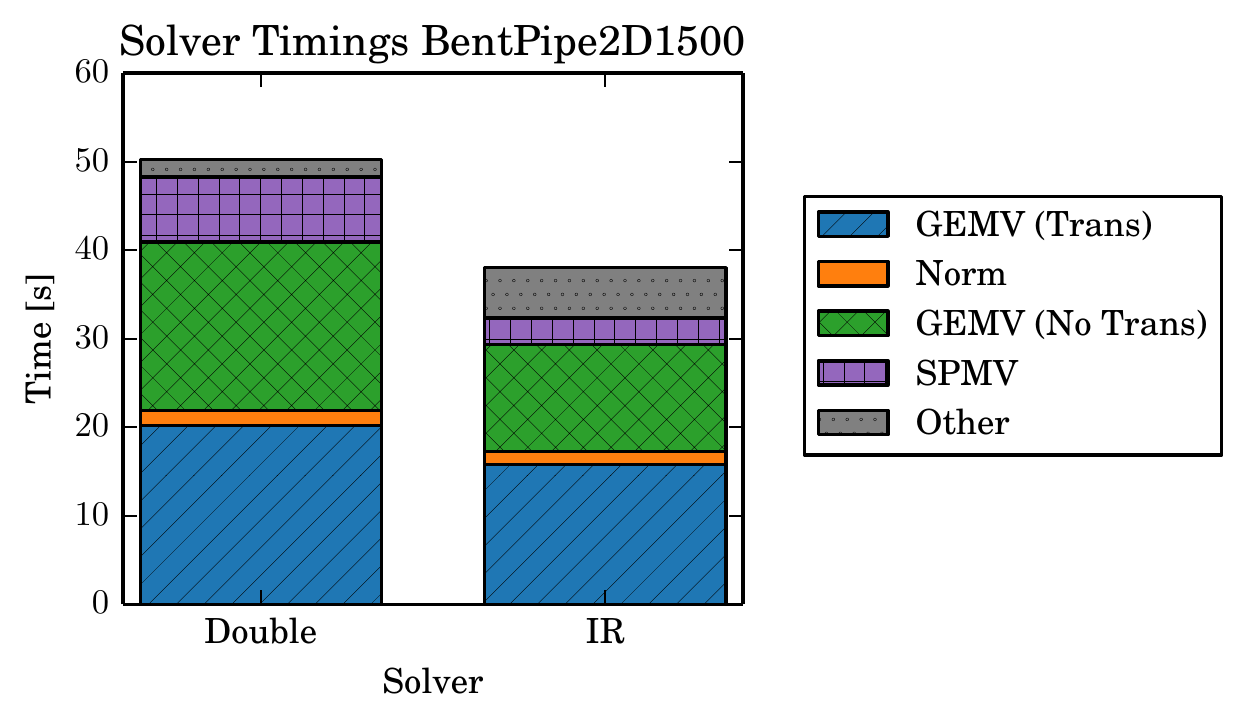}
  \caption{Solve times for GMRES$(50)$ double (left) and IR (right) for the matrix BentPipe2D1500. Each bar represents total solve time, split up to give a breakdown of time spent in different kernels. The ``Other" portion represents timing for small dense (non-GPU) operations and, for GMRES-IR, computing residuals in fp64.} %TODO: See if we can add a bar split in the IR for time to compute double residuals. 
  \label{fig:BentPipeTiming}
\end{figure}

Figure \ref{fig:BentPipeTiming} and Table \ref{tab:BentPipeSpeedup} show the solve times and speedup of the GMRES double and IR solvers, split over different kernels. Solve times do not include time required to copy the matrix $A$ from fp64 to fp32 %\JL{Should I capitalize the fp?  Or italicize or something?} 
at the beginning of GMRES-IR. By this measure, GMRES-IR gives $1.32\times$ speedup over the solve time of GMRES double. The two GEMV kernels give $1.28$ to $1.57\times$ speedup, but the SpMV gives a spectacular $2.48\times$ speedup! %\CG{I think you checked this already: we are not mixing SpMV implementations and using KK for float but CuSPARSE for double? } No- it is using KK for all calls, double and float. 
%Even so, %since the SpMV kernel only comprised $15\%$ of the original solve time, the contribution to overall speedup is small. 
%We will discuss the SpMV speedup further in Section \ref{sec:SpMV}. %This is said again later.
The bar segment in Figure \ref{fig:BentPipeTiming} labeled ``other" indicates time solving the least squares problems and performing other non-GPU operations in GMRES. For GMRES-IR, it also includes computation of the new residual in double precision. %TODO: separate that part out in the graph. 

%TODO Edit this table to match bar graph
%TODO Do we really need both the graph and table?
% Table generated by Excel2LaTeX from sheet 'BentPipePaper'
\begin{table}[htbp]
 \centering
 \caption{Speedup of different kernels for the matrix BentPipe2D1500.}
  \begin{tabular}{lrrr}
   & \multicolumn{1}{l}{Double Belos} & IR Belos & Speedup \\
  %& \multicolumn{1}{l}{Double Belos} & \multicolumn{1}{l}{IR Belos} & \\
  %   &    &    & \multicolumn{1}{l}{Speedup} \\
  
  %\textbf{Total solve time for convergence: } & \textbf{50.23} & \textbf{34.05} & 1.47518355 \\
  \hline
  \textbf{GEMV (Trans)} & 20.20 & 15.78 & 1.28 \\
  \textbf{Norm} & 1.72 & 1.49 & 1.15 \\
  \textbf{GEMV (no Trans) } & 19.01 & 12.10 & 1.57\\
  \textbf{Total Orthogonalization } & 41.85 & 30.30 & 1.38 \\
  \textbf{SpMV} & 7.33 & 2.95 & 2.48 \\
  \textbf{Total Time} & 50.26 & 38.03 & 1.32\\
  \end{tabular}%
  
  %\CG{Maybe list time/call since the number of calls isn't necessarily identical? }
 \label{tab:BentPipeSpeedup}%
\end{table}%

\begin{figure}
  \centering
  \includegraphics[width=0.38\textwidth,keepaspectratio]{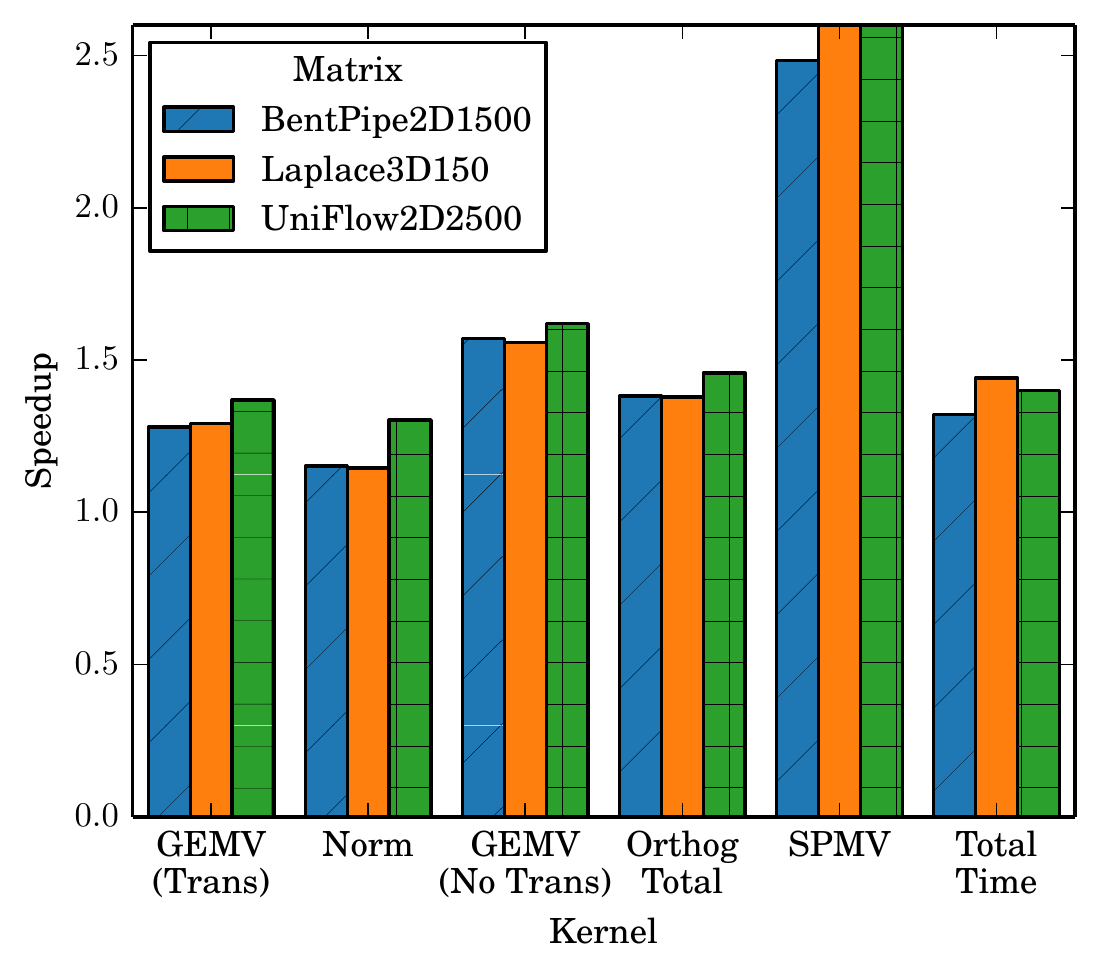}
  \caption{Speedups for different kernels going from GMRES double to GMRES-IR over three different PDEs. (Note that this is speedup of the total time spent in each kernel in GMRES double vs GMRES-IR. This is not a per-call comparison.)}
  \label{fig:KernelSpeedup3Mats}
\end{figure}

In Figure \ref{fig:KernelSpeedup3Mats}, we graph kernel speedups for the previous problem and two additional matrices: a 3D Laplacian and the matrix UniFlow2D2500 from Section \ref{subsec:IRvsFD}. %The latter is another convection-diffusion problem from the Galeri package with $n=6.25$ million. %Now we've already introduced this. 
(See Table \ref{tab:largeTestSet} for additional problem statistics.)
Note that these bars show the speedup of the entire time GMRES double spends in a kernel over the entire time GMRES-IR spends in the same kernel. Since GMRES-IR needs a few extra iterations (and kernel calls) beyond what GMRES double needs to converge, this is not a per-call time comparison. Even so, speedups for a per-call comparison are very similar to those presented in Figure \ref{fig:KernelSpeedup3Mats}.
It is interesting to note that the kernel speedups are relatively consistent across the three problems. In particular, the SpMV kernel improves by $2.4$ to $2.6$ times in all three cases. This occurs due to near-perfect L2 cache reuse for the right-hand side vector with SpMV float, while there is a high L2 cache miss rate for SpMV double. We will discuss SpMV speedup further in Section \ref{sec:SpMV}.
The total solve times to convergence for the three problems improve by $24$ to $36\%$. 
%\CG{[The speedup from table 1 is $1.32\times$, fig 5 looks similar, but $1/1.32=0.75$, so not sure 40\% improvement is correct.]}

%\CG{Is there an explanation for the GEMV difference for transpose vs no-transpose? Or are these simply calls with different matrix/vector sizes?} \JL{They are called with different sizes, but they are also completely different kernels. Transpose is Seher's dot-based GEMM. Else is CuBlas. }

\subsection{Convergence and Kernel Speedup for Preconditioned GMRES vs GMRES-IR}
\label{subsec:precon_compar}
In this section we compare three preconditioning options.
%: GMRES-double with double precision preconditioning, GMRES-double with single precision preconditioning, and GMRES-IR with single precision preconditioning.
The matrix $A$ is a 2D Laplacian over a stretched grid with $nx = 1500$. %Thus, we have $n=2.25$ million as with the BentPipe2D1500 matrix, but here the bandwidth of the matrix is larger. %True, but how does that make a difference? This prob is also symmetric. 
It has a large condition number, so GMRES$(50)$ cannot converge without preconditioning. 
%\JL{TODO Maybe use deg 20 instead of 40? Nicer gap between Single prec and IR....}
%Having only the preconditioner in single precision gives a new method to mix precisions. 
We apply a GMRES polynomial preconditioner \cite{LoePPTrilinos} of degree $40$, using a) GMRES-fp64 with fp64 preconditioning, b) GMRES-fp64 with fp32 preconditioning, and c) GMRES-IR with fp32 preconditioning. Here ``fp32 preconditioning" indicates that the polynomial is both computed and applied in single precision. 

\begin{figure}
  \centering
  \includegraphics[width=0.45\textwidth,keepaspectratio]{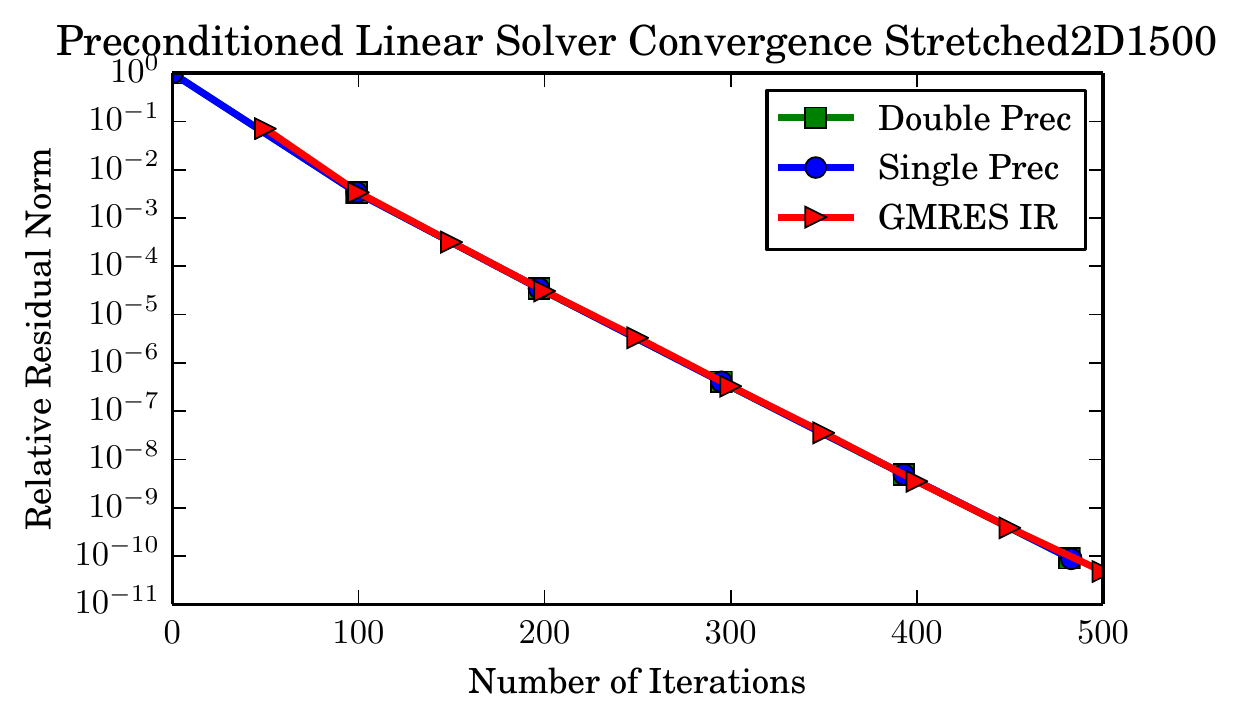}
  \caption{Convergence of the Stretched2D1500 problem with a degree $40$ polynomial preconditioner. Squares indicate fp64 preconditioning, circles fp32 preconditioning, and triangles GMRES-IR with fp32 preconditioning.}
  \label{fig:polyConv}
\end{figure}

\begin{figure}
  \centering
  \includegraphics[width=0.35\textwidth,keepaspectratio]{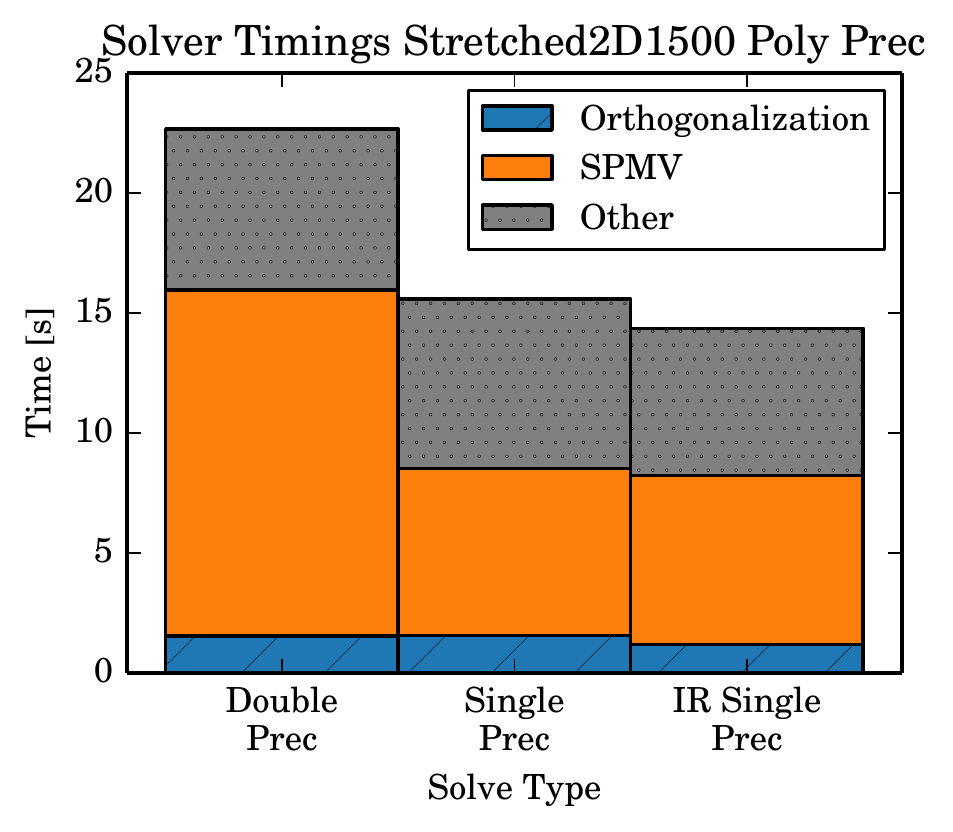}
  \caption{Solve times for polynomial preconditioned GMRES using polynomial degree $40$. The bar on the left shows solve time for fp64 GMRES, the bar in the middle shows fp64 GMRES with an fp32 polynomial, and the bar on the right gives timings for GMRES-IR with fp32 polynomial preconditioning.}
  \label{fig:PolyPrecKernelSpeedup}
\end{figure}

Figure \ref{fig:polyConv} demonstrates that, just as before, the problems with fp32 preconditioning converge very similarly to GMRES in all fp64. 
%In fact, GMRES preconditioned in single precision converges within one iteration of the all-double precision solver. 
Figure \ref{fig:PolyPrecKernelSpeedup} shows solve times for all three configurations. Times do not include creation of the polynomial preconditioner, which was $0.5$ seconds or less for all cases. 
Similar to Figure \ref{fig:BentPipeTiming}, the ``other" portion of each bar indicates time spent in dense matrix operations, vector additions for the polynomial, and computation of double-precision residuals in GMRES-IR. Since the SpMV constitutes the majority of kernel calls in the polynomial preconditioner and gets large speedup, the total SpMV time drops significantly in single precision as opposed to double. Time spent in ``other" operations, however, increases slightly due to the casting operations required to multiply an fp32 matrix polynomial with an fp64 vector. %However, with degree $60$, the single precision preconditioned GMRES needs $42$ additional iterations over the all-double precision GMRES. 
Ultimately, GMRES-IR gives $1.58\times$ speedup over GMRES double. Even when testing other polynomial degrees, the fp32 preconditioned GMRES gives reasonable speedup over the all-double precision GMRES, but run times are never faster than those of GMRES-IR.

Unlike previous examples where solve time was dominated by orthogonalization, polynomial preconditioning shifts the cost toward the sparse matrix-vector product. Here, the SpMV gets about $2\times$ speedup going from fp64 to fp32. %(This is likely less than the $2.5$ times speedup of earlier since the Stretched2D matrix has larger bandwidth \CG{[number of nonzeros per row, not bandwidth]} than the examples in Figure \ref{fig:KernelSpeedup3Mats}.) %Is that true???
Note that in the previous example (Figure \ref{fig:BentPipeTiming}), the BentPipe SpMV kernel only comprises $15\%$ of the fp64 solve time, so the $2.5\times$ SpMV speedup only removes $4.4$ seconds from the original solve time of $50$ seconds. In this stretched Laplacian problem, the SpMV comprises $64\%$ of the total solve time for fp64, so the improvement in SpMV time provides $32\%$ of the ultimate speedup in GMRES-IR. Polynomial preconditioning allows us to take advantage of the large speedup from applying the SpMV in lower precision. 
%Here, with the degree $10$ polynomial, the SpMV comprises $y\%$ of the double precision solve time, and contributes uuu reduction in the final kernel. For the degree $60$ polynomial, this is even more significant. 

% \JL{TODO: Where to move this?? Or maybe just cut it.} \sr{Keep it if space is fine.}
% We also note that the degree $60$ problem has a slightly lower speedup than the other polynomial degrees since it takes $53$ more iterations than the double precision problem. In the future, it might be helpful to stop the iterative refinement solver and check for convergence sooner than the next multiple of $50$ iterations. The solver had almost converged to $1\mathrm{e}{-}10$ at the $6$th restart with $300$ iterations, but then had to complete another cycle of $50$ iterations which were not entirely necessary before recomputing the residuals. 

While this analysis has only covered polynomial preconditioning, we believe that the following concepts will also extend to many other preconditioners: a) Convergence of problems preconditioned in fp32 will typically follow convergence of fp64 preconditioning; b) Using an fp32 preconditioner with fp64 GMRES will typically improve solve time over using the same preconditioner in fp64, but perhaps not as much as applying that preconditioner within GMRES-IR; and c) Preconditioning allows users to take advantage of kernels that have large speedup in lower precisions. 

\subsection{Matrix Structure, Cache Reuse, and SpMV Performance}
\label{sec:SpMV}

The roughly $2.5\times$ speedup of the sparse matrix-vector product in the previous examples requires deeper explanation. Intuitively, one might expect that changing the working precision from fp64 to fp32 should give at most 1.5 to two times speedup since we are reducing the memory requirement by almost half. We assume the integer index type stays the same. If we halve the floating point data size and the index size stays the same, %(typically 64 bits, like double), 
then one might expect at most $1.5 \times$ speedup. Below we explain how lower precision can improve cache reuse and give greater than $1.5$ or even $2$ times speedup. 

Note that the SpMV kernel called in all previous examples is an implementation native to Kokkos Kernels; we do not employ CuSparse for SpMV (though CuBlas may be called in other operations). The SpMV kernel is memory-bound, so the limiting factor in speed is how fast data can be moved through the memory hierarchy. 
Recall that storing a double requires $8$ bytes of memory and that both integers and floats require $4$ bytes of memory. Each of our matrices is stored in Compressed Sparse Row (CSR) format. With NVIDIA profiling tools, we observed that the L2 cache hit rate for the float SpMV was almost twice the hit rate for the double SpMV. This appears to be due to ``perfect caching" of the right-hand side vector $x$. Below we give a calculation to explain how this caching effect can account for $2.5\times$ speedup. 

Suppose that $A$ has $w$ nonzero elements per row and $n$ rows (so $nnz=w*n$) and that we are computing $Ax=y$. With the CSR matrix storage format, we have two vectors of length $nnz$ [one for the values of $A$ (denoted $A_{val}$) and another for the column indices (denoted $colId$)] and a vector of row pointers of length $n+1$. For this calculation, we ignore reads of the vector of row pointers and writes to $y$ since they account for only a small fraction of all memory traffic. To compute the dot product for each element in the solution vector $y$, we have to read one row of nonzeros $A$ and $w$ elements of $x$ which correspond to their locations. Thus the first dot product is 
\[\sum_{i=0}^{w-1} A_{val}[i] * x[colId[i]]. \]
Suppose now that in fp64, there is no cache reuse for the $x$ vector; we have to reread each element from device memory to cache every time we need it. Then to compute the SpMV, for each nonzero element in $A$ we read one double from $A$, one int (for $colId[i]$), and another double from $x$. In that case, the total number of reads to cache is 
\[n*w*[size(int) + 2* size(double)] = 20wn.\]
Next, we suppose that in fp32 there is ``perfect caching" of the $x$ vector. In other words, we only have to read $x$ from device memory once, and after an element is read into cache, it stays there until we do not need it any longer. In that case, the total number of reads to cache is 
\[n*w*[size(int) + size(float)] +n*size(float) = (8w+4)n.\]
Then the speedup going from double to float is
\[ \frac{20wn}{(8w+4)n} = \frac{5w}{2w+1}.\]
This ratio quickly approaches $2.5$ as $w$ grows. For the matrices in Section \ref{subsec:unprecon_compar}, the speedup as predicted by the model is slightly lower. Matrices UniFlow2D2500 and BentPipe2D1500 have $5$ nonzeros per row, so the expected speedup is $2.27\times$. With the Laplace3D150 matrix that has $7$ nonzeros per row, the expected speedup from this model is $2.33\times$. The observed speedup in all three cases was slightly higher than expected, probably due to additional improvements in $L1$ cache use. %TODO: Is that last "probably" statement true?

Additional experiments have confirmed this model: for nicely structured matrices, $2.5\times$ speedup can result from perfect cache reuse for $x$ in fp32, while some $x$ vector elements must be re-read into cache for fp64. 
Note that if $A$ has larger bandwidth, elements of $x$ may be accessed with less spatial locality, so $2.5\times$ speedup is not expected. For additional study of SpMV in multiple precisions, see \cite{AhmadMPSpMV}.

%Experiments have suggested that the matrix structure (e.g. number of non-zeros per row) can have a large impact on the speedup of SpMV operations going from double precision to single precision. . Some matrices get over 2x speedup due to cache reuse, and other problems get almost no speedup at all. (The paper \cite{AhmadMPSpMV} had similar variations in speedup with their algorithms.) Since the SpMV is such a crucial kernel in GMRES (and in polynomial preconditioning) its contribution to overall speedup is significant. 

%\CG{Might be a good idea to run some 1D example. There the access pattern should be super nice, so we should get close to 2x speedup?}

\subsection{Choosing a Restart Size for GMRES-IR}
\label{subsec:restart}
Here we demonstrate an interesting case where choosing a small restart size for GMRES gives improved performance (in double and mixed precisions) over a large subspace size.
The authors of \cite{LindquistGMRES} devote many experiments to determining the best restart strategy for GMRES-IR. Their strategy is to pick the restart size that allows the inner low precision GMRES to converge as far as possible before restarting. This means picking the largest subspace possible before convergence in the inner solver stalls. Here we demonstrate a further example with matrix BentPipe2D1500 where the large matrix size causes orthogonalization costs to dominate the solve time. Thus, a smaller restart size is more beneficial for this problem. 

We test a variety of restart sizes; Table \ref{tab:ManySubsBentPipe} gives the solve times and iteration counts. In each case, GMRES-IR still gives speedup of $1.20\times$ to $1.40\times$ over GMRES double. 
% Table generated by Excel2LaTeX from sheet 'BentPipePaper'
\begin{table}[htbp]
 \centering
 \caption{BentPipe2D1500 Convergence for Many Restart Sizes}
  \begin{tabular}{llrrrr}
  \multicolumn{1}{p{4.085em}}{\textbf{Subsp }} & \multicolumn{2}{c}{\textbf{GMRES Double}} & \multicolumn{2}{c}{\textbf{GMRES-IR}} & \\
  \multicolumn{1}{p{4.085em}}{\textbf{Size}} & Iters & Solve Time & Iters & Solve Time & \multicolumn{1}{l}{Speedup} \\
  \hline
  25  & 13795 & 38.63 & 13925 & 31.74 & 1.22 \\
  50  & 12967 & 50.26 & 13150 & 38.03 & 1.32 \\
  100  & 12009 & 74.24 & 12100 & 51.88 & 1.43 \\
  150  & 11250 & 95.82 & 12450 & 72.01 & 1.33 \\
  200  & 10867 & 117.80 & 12400 & 90.77 & 1.30 \\
  300  & 10491 & 164.60 & 12600 & 133.60 & 1.23 \\
  400  & 10274 & 209.80 & 12400 & 174.10 & 1.21 \\
  
  \end{tabular}%
 \label{tab:ManySubsBentPipe}%
 
 %\JL{Remove comment about subsp 100 IR taking fewer iters.}
\end{table}%
Although the iteration count for the fp64 solver decreases as the subspace gets larger, the solve time increases. The large subspace size causes orthogonalization costs to increase and dominate the solve time more and more. Figure \ref{fig:BentPipeTiming} demonstrates the proportion of total orthogonalization costs (GEMV Trans $+$ Norm $+$ GEMV no Trans) for restart size $50$. In that figure, orthogonalization consumes $83\%$ of solve time for the fp64 solver and $80\%$ of solve time for GMRES-IR. As the restart size increases, the proportion of solve time for SpMVs and non-orthogonalization operations gets squeezed out. With a restart size of $400$, orthogonalization takes $97\%$ of solve time for both the double precision and IR solvers. 

The smallest restart size of $25$ also gives the best solve time for GMRES-IR. Like the double precision solver, GMRES-IR benefits from reduced orthogonalization costs with the small restart length. Observe that, contrary to the strategy in \cite{LindquistGMRES}, size $25$ gives us the fastest solve time even though the single precision inner solver convergence is not near stalling. Even for the largest restart size of $400$, the residuals of the inner solver do not appear to have stalled; they are all on the order of $0.1$. Typically an fp32 solver can converge to near $10^{-5}$ without any fp64 refinement. %In all cases, the GMRES-IR solver is faster than the double precision solver restarted at the same interval. 
%Note, however, that the $1.44\times$ speedup for subspace size $100$ may be slightly inflated: In this case, GMRES-IR actually takes fewer iterations than GMRES double, meaning that some of the speedup comes from simply performing fewer operations rather than from using float precision. 

As the inner fp32 GMRES solver is restarted (and refined) less frequently, the gap between the iteration count needed for GMRES-IR and GMRES double convergence widens. Note that, while unusual, accumulated rounding errors can occasionally help GMRES-IR to need fewer iterations to converge than GMRES double. This phenomenon did not happen in the presented median-time runs, but it was occasionally present in other runs of the experiment. Nevertheless, the GMRES-IR solver consistently gives performance improvement over the GMRES double solver. 

\begin{figure}
  \centering
  \includegraphics[width=0.5\textwidth,keepaspectratio]{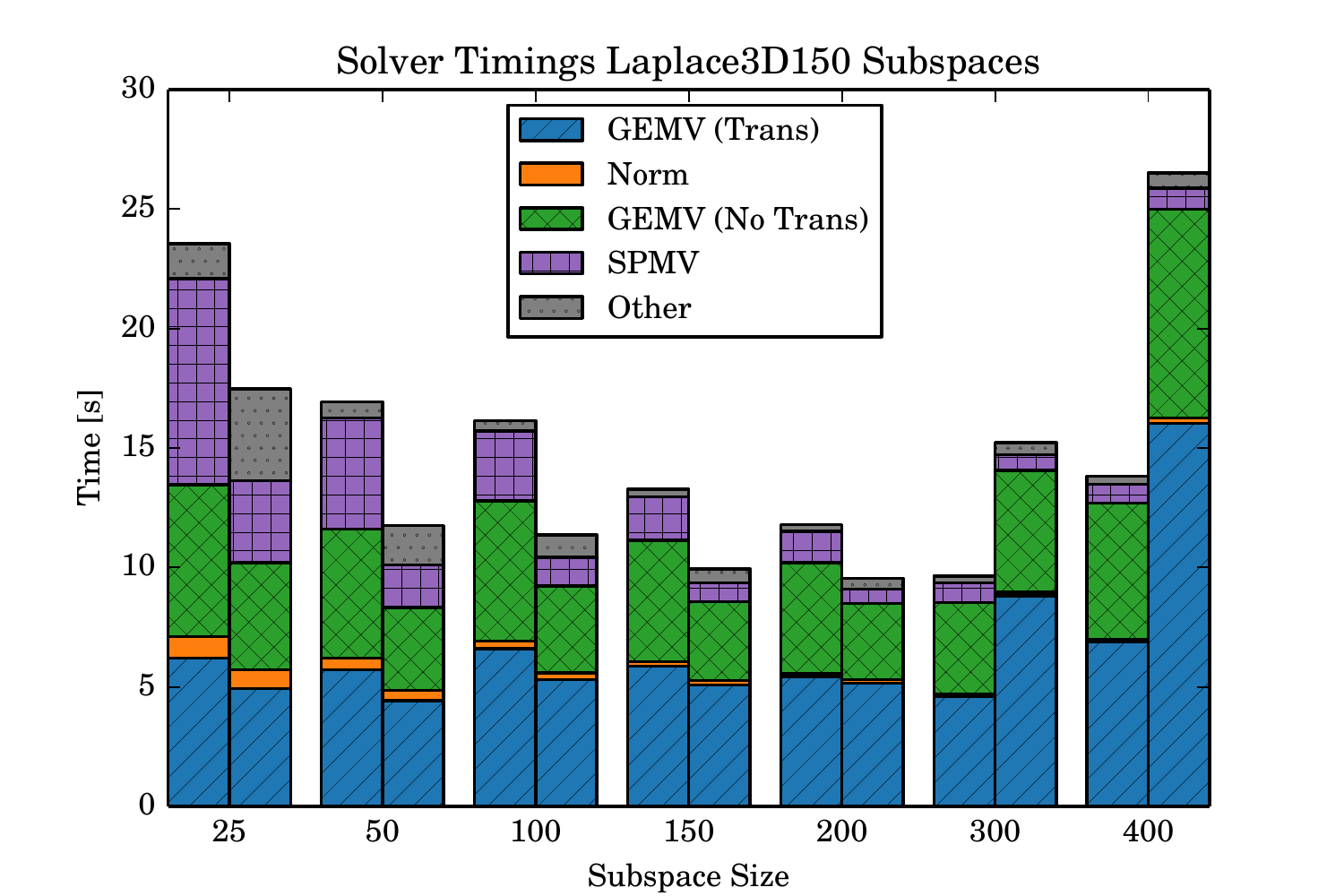}
  \caption{Total solve times for different GMRES restart lengths for the matrix Laplace3D150. For each restart size, the left bar indicates solve time for GMRES (double) and the right bar gives solve time for GMRES-IR. }
  \label{fig:LaplSubspMulti}
\end{figure}

Next, we show an example using a 3D Laplacian where GMRES-IR does not give speedups at large subspace sizes. %The matrix is a 3D Laplacian with $150$ grid points in each direction, so $n=3{,}375{,}000$ and $nnz = 23{,}490{,}000$. %Table \ref{tab:LaplSubspMulti} gives the solve times and iteration count for different GMRES restart sizes. 
Results are in Figure \ref{fig:LaplSubspMulti}, where bars indicating solve time are broken down into the times for particular kernels. 
For restart sizes up to $200$, the GMRES-IR solver gives $19\%$ to $31\%$ improvement in solve time over GMRES double. % Formula for percent speedup: percent = 1-(fastTime/slowTime). 
However, with larger subspaces, the iterative refinement solver needs so many additional iterations over GMRES double that we do not see any speedup. With size $300$, GMRES double needs $433$ iterations compared to $900$ iterations for GMRES-IR. For subspace size $400$, GMRES-IR needs almost three times as many iterations as GMRES double. In the experiments for both of these large subspace sizes, we see strong evidence of stalled convergence in the single precision solver; several residuals are on the order of $10^{-7}$. Slowdown comes with GMRES-IR because the double precision residual is updated so infrequently; the inner solver is taking extra iterations without making progress towards the solution. 

Ultimately, the fastest solve time is with GMRES-IR and a subspace size of $200$ (though the timing of GMRES$(300)$ in double was faster on some runs). %Beginning with subspace size $400$, the increased orthogonalization costs begin to outweigh the small decrease in iterations. 
It should be noted that for larger versions of this PDE matrix, attempting to use a subspace size of $300$ results in an out-of-memory error on the GPU. Thus, GMRES-IR likely gives the most practical gains in terms of solve time for large problems.  

% Moving this table earlier in text so it doesn't appear after references.
% Table generated by Excel2LaTeX from sheet 'Matrix list'
\begin{table*}[htbp]
 \centering
   \caption{Timings and iteration counts for GMRES double and GMRES-IR for a variety of Suitesparse and Galeri matrices.}
  \begin{tabular}{rlrrll|rr|rr|c}
     &    &    &    &    &    & \multicolumn{2}{c}{\textbf{Double}} & \multicolumn{2}{c}{\textbf{IR}} & \\
  \multicolumn{1}{l}{\textbf{UF ID}} & \textbf{Matrix Name} & \multicolumn{1}{l}{\textbf{N}} & \multicolumn{1}{l}{\textbf{NNZ}} & \textbf{Symm} & \textbf{Prec} & \multicolumn{1}{l}{\textbf{Time}} & \multicolumn{1}{l}{\textbf{Iters}} & \multicolumn{1}{l}{\textbf{Time}} & \multicolumn{1}{l}{\textbf{Iters}} & \multicolumn{1}{l}{\textbf{Speedup}} \\ \hline
  2266 & atmosdmodj & 1,270,432 & 8,814,880 & n   &    & 5.12 & 1740 & 3.78 & 1750 & 1.35 \\
  1849 & Dubcova3 & 146,698 & 3,636,643 & spd  &    & 1.15 & 1131 & 1.05 & 1150 & 1.10 \\
  895  & stomach & 213,360 & 3,021,648 & n   &    & 0.51 & 359  & 0.52 & 400  & 0.98 \\
  1367 & SiO2 & 155,331 & 11,283,503 & y   &    & 18.23 & 17385 & 16.86 & 17600 & 1.08 \\
  1853 & parabolic\_fem & 525,825 & 3,674,625 & spd  &    & 41.77 & 27493 & 45.34 & 36600 & 0.92 \\
  894  & lung2 & 109,460 & 492,564 & n   & J 1  & 0.46 & 206  & 0.49 & 250  & 0.94 \\
  1266 & hood & 220,542 & 9,895,422 & spd  & J 42 & 13.98 & 5762 & 9.04 & 5000 & 1.55 \\
  805  & cfd2 & 123,440 & 3,085,406 & spd  & p 25 & 6.05 & 1092 & 4.55 & 1100 & 1.33 \\
  2649 & Transport & 1,602,111 & 23,487,281 & n   & p 25 & 8.35 & 339  & 8.73 & 450  & 0.96 \\
  1431 & filter3D & 106,437 & 2,707,179 & y   & p 25 & 25.24 & 4449 & 18.12 & 4450 & 1.39 \\
     & BentPipe2D1500 & 2,250,000 & 11,244,000 & n   &    & 50.26 & 12967 & 38.03 & 13150 & 1.32 \\
     & UniFlow2D2500 & 6,250,000 & 31,240,000 & n   &    & 29.62 & 2905 & 21.17 & 3000 & 1.40 \\
     & Laplace3D150 & 3,375,000 & 23,490,000 & spd  &    & 16.93 & 2387 & 11.75 & 2400 & 1.44 \\
     & Stretched2D1500 & 2,250,000 & 20,232,004 & spd  & p 40 & 22.66 & 482  & 14.37 & 500  & 1.58 \\
     \hline
  \end{tabular}%
  
  \vspace{0.05in}
     \textnormal{In the ``Symm" column, 'y' indicates a symmetric matrix and 'spd' indicates a symmetric positive definite matrix. In the ``Prec" column, ``J $k$" indicates a block Jacobi preconditioner with block size $k$, and ``p $k$" indicates a polynomial preconditioner of degree $k$.} 
 \label{tab:largeTestSet}%
\end{table*}%

\subsection{Choosing Preconditioner Complexity for Multiprecision}
\label{subsec:poly}

Here a new example shows that while using a single precision preconditioner for double precision GMRES can often work well, users need to keep a caveat in mind: The more computationally intensive the fp32 preconditioner is, the more opportunities there will be for significant rounding errors to accumulate. Furthermore, preconditioners that work well in fp64 may be unstable in fp32. We demonstrate this by polynomial preconditioning a $3$D Laplacian that has $200$ grid points in each direction. We test polynomial degrees that are multiples of $10$ up to $70$. When we apply the polynomial and all other operations in fp64, the GMRES solver always converges successfully. Then we apply the polynomial preconditioner in fp32 and perform all other GMRES calculations in fp64. For the degree $10$ polynomial the solver converges, just as it does in all double precision. However, for higher degree polynomials, the implicit residual (that which results from applying Givens rotations to the matrix $H$ from the Arnoldi relation) diverges from the explicit residual (computed by forming $\hat{x}$ and calculating $\|b-A\hat{x}\|_2$). In the Belos solvers library, divergence of the implicit and explicit residuals is denoted as a ``loss of accuracy" of the solver. In essence, the solver gives a ``false positive" signal of convergence. 
Thus, when applying the polynomial becomes more computationally expensive (with high degrees), accumulated rounding errors prevent the solver from converging with respect to the true residual $\|b-A\hat{x}\|_2$. It is also likely that this preconditioner becomes ill-conditioned far more quickly in single precision than in double. Scientists should bear these possibilities in mind when applying any fp32 preconditioner to fp64 GMRES. GMRES-IR is less likely to suffer from diverging implicit and explicit residuals since it performs a correction with the true residual at each restart. Forcing such a correction in GMRES double after the loss of accuracy is likely to fix the problem with the fp32 preconditioner. Applying this fix is possible in Belos, but it is not yet built in.

\subsection{Large Test Set from SuiteSparse}
\label{subsec:suitesparse}

Finally, we validate the prior analysis with additional examples. Several matrices from the SuiteSparse matrix collection \cite{DavisSparseCollect} are tested with GMRES double and GMRES-IR. Results are in Table \ref{tab:largeTestSet}. The first five matrices do not have preconditioning. The next two matrices are reordered with a reverse Cuthill-McKee ordering before applying block Jacobi preconditioners with block sizes of $1$ and $42$, respectively. The next three matrices in the table use polynomial preconditioners of degree $25$. At the end of the table, we repeat the earlier results of Section \ref{sec:experiments} for completeness. 

Results for speedup from GMRES-IR are mixed. Broadly generalizing, GMRES-IR seems more likely to give speedup for matrices which need many hundreds or thousands of iterations to converge. For matrices which need very few iterations to converge with double precision GMRES, the cost of the additional iterations required by GMRES-IR seems to outweigh any gain from incorporating single precision. The \textit{parabolic\_fem} problem needs further investigation; unlike in Section \ref{subsec:unprecon_compar}, GMRES-IR convergence for this matrix quickly diverges from the convergence of GMRES double. In problems where we do see speedup, the values vary from $1.08\times$ to $1.58\times$. Typically GMRES$(50)$-IR needs a few more iterations to converge than GMRES$(50)$ double, but the \textit{hood} matrix is a counterexample. For the \textit{hood} matrix, roundoff errors allow GMRES-IR to converge with $762$ fewer iterations than GMRES double, giving us a higher speedup than expected from simply switching to a lower working precision. %For the \textit{Transport} problem, unfortunately, GMRES-IR needs more than two times as many iterations as GMRES double, and this results in overall slowdown.

\section{Conclusions and Future Work}
In this work, we evaluated two different approaches for multiprecision GMRES. We found GMRES-IR to be the best choice. GMRES-IR is a flexible algorithm for incorporating lower precision calculations into GMRES while maintaining double precision accuracy of the final solution. Its convergence is similar to double precision GMRES and it often provides speedup for large problems which require many iterations. We observed similar results with preconditioned GMRES-IR as well. We analyzed the speedup at the individual kernel level and recommended guiding principles for selecting solver parameters. In the future, we plan to make GMRES-IR available for Trilinos users in the Belos solvers package. We believe this could replace standard (all double) GMRES in many applications. We will further study its capabilities on multiple GPUs and with other preconditioners. Since Kokkos is enabling support for half precision, we will also study ways to incorporate a third level of precision into the GMRES-IR solver while maintaining high accuracy. This may further improve solve time for applications problems.

\section*{Acknowledgment}
Thanks to Christian Trott and Luc Berger-Vergiat for helping develop the model in Section V-D. We also thank the referees for their may useful suggestions. 
This research was supported by the Exascale Computing Project (17-SC-20-SC), a collaborative effort of the U.S. Department of Energy Office of Science and the National Nuclear Security Administration.

%\section*{References}

%TODO- How to handle works like Trilinos and MP overview with many many authors?
\bibliographystyle{IEEEtran}
% I created a copy of MultiPrecision.bib with all URLs, DOIs etc removed.
% Just change back to MultiPrecision.bib if we want them again.
\bibliography{IEEEabrv,MultiPrecisionNoURL.bib}

\end{document}